\RequirePackage{ifpdf}
\ifpdf 
\documentclass[pdftex]{sigma}
\else
\documentclass{sigma}
\fi

\begin{document}

\allowdisplaybreaks

\renewcommand{\thefootnote}{$\star$}

\renewcommand{\PaperNumber}{019}

\FirstPageHeading

\ShortArticleName{Besov-Type Spaces on $\mathbb{R}^d $ and Integrability for the Dunkl Transform}

\ArticleName{Besov-Type Spaces on $\mathbb{R}^d $ and Integrability\\ for the Dunkl Transform\footnote{This paper is a contribution to the Special
Issue on Dunkl Operators and Related Topics. The full collection
is available at
\href{http://www.emis.de/journals/SIGMA/Dunkl_operators.html}{http://www.emis.de/journals/SIGMA/Dunkl\_{}operators.html}}}

\Author{Chokri ABDELKEFI~$^\dag$, Jean-Philippe ANKER~$^\ddag$,
Feriel SASSI~$^\dag$ and Mohamed SIFI~$^\S$}

\AuthorNameForHeading{C. Abdelkef\/i, J.-Ph. Anker, F. Sassi and M. Sif\/i}

\Address{$^\dag$~Department of Mathematics, Preparatory
Institute of Engineer Studies of Tunis,\\
\hphantom{$^\dag$}~1089 Monf\/leury Tunis, Tunisia}
\EmailD{\href{mailto:chokri.abdelkefi@ipeit.rnu.tn}{chokri.abdelkefi@ipeit.rnu.tn},
\href{mailto:feriel.sassi@ipeit.rnu.tn}{feriel.sassi@ipeit.rnu.tn}}

\Address{$^\ddag$~Department of Mathematics, University of Orleans
$\&$ CNRS, Federation Denis Poisson\\
\hphantom{$^\ddag$}~(FR 2964), Laboratoire MAPMO
(UMR 6628), B.P. 6759, 45067 Orleans cedex 2, France}
\EmailD{\href{mailto:Jean-Philippe.Anker@univ-orleans.fr}{Jean-Philippe.Anker@univ-orleans.fr}}

\Address{$^\S$~Department of Mathematics, Faculty of Sciences
of Tunis, 1060 Tunis, Tunisia}
\EmailD{\href{mailto:mohamed.sifi@fst.rnu.tn}{mohamed.sifi@fst.rnu.tn}}

\ArticleDates{Received August 28, 2008, in f\/inal form February 05,
2009; Published online February 16, 2009}

\Abstract{In this paper, we show the inclusion and the density of
the Schwartz space in Besov--Dunkl spaces and we prove an
interpolation formula for these spaces by the real method. We give
another characterization for these spaces by convolution. Finally,
we establish further results concerning integrability of the Dunkl
transform of function in a~suitable Besov--Dunkl space.}


\Classification{42B10; 46E30; 44A35}

\section{Introduction}

We consider the dif\/ferential-dif\/ference operators $T_i$, $1 \leq i
\leq d$, on $\mathbb{R}^d$, associated with a positive root system
$R_+$ and a non negative multiplicity function $k$, introduced by
C.F.~Dunkl in~\cite{[9]} and called Dunkl operators (see next section).
These operators can be regarded as a generalization of partial
derivatives and lead to generalizations of various analytic
structure, like the exponential function, the Fourier transform, the
translation operators and the convolution (see~\cite{[8], [10], [11], [16], [17],
[18], [19], [22]}). The Dunkl kernel $E_k$ has been introduced by C.F.~Dunkl in~\cite{[10]}. This kernel is used to def\/ine the Dunkl transform
$\mathcal{F}_k$. K. Trim\`eche has introduced in~\cite{[23]} the Dunkl
translation operators~$\tau_x$,
 $x\in \mathbb{R}^d$, on the space of inf\/initely
dif\/ferentiable functions on $\mathbb{R}^d$. At the moment an
explicit formula for the Dunkl translation operator of function
$\tau_x(f)$ is
 unknown in general. However, such formula is known when $f$ is a radial function and the $L^p$-boundedness of $\tau_x$
 for radial functions is established. As a result, we have
 the Dunkl convolution $\ast_k$.

There are many ways to def\/ine the Besov spaces (see~\cite{[6],[15], [21]}) and the
Besov spaces for the Dunkl operators (see~\cite{[1], [2], [3], [4], [14]}).
 Let $ \beta>0$, $1\leq p , q \leq + \infty$, the Besov--Dunkl space
denoted by $ \mathcal{B}\mathcal{D}_{p,q}^{\beta,k} $ in this paper,
is the subspace of functions $f \in L^p_k(\mathbb{R}^d)$ satisfying
\[
\|f\|_{\mathcal{B}\mathcal{D}_{p,q}^{\beta,k}}=\bigg(\sum_{j\in\mathbb{Z}}(2^{j\beta}\|\varphi_{j
}\ast_k f\|_{p,k})^q \bigg)^{\frac{1}{q}}<+\infty  \qquad \mbox{if}\quad
q < +\infty
\]
 and
\[
\|f\|_{\mathcal{B}\mathcal{D}_{p,\infty}^{\beta,k}}= \sup_{j\in
\mathbb{Z}}2^{j\beta}\|\varphi_{j }\ast_k
f\|_{p,k}<+\infty\qquad  \mbox{if} \quad q = +\infty ,
\] where $
(\varphi_{j})_{j\in \mathbb{Z}}$ is a sequence of functions in $
\mathcal{S}(\mathbb{R}^d)^{\rm rad}$ such that
\begin{itemize}\itemsep=0pt
\item[$(i)$] $\mbox{supp}\,\mathcal{F}_k(\varphi_{j}) \subset A_j=\big\{x\in
\mathbb{R}^d\,;\,2^{j-1}\leq \|x\|\leq 2^{j+1}\big\}$ for $j\in
\mathbb{Z}$;
\item[$(ii)$] $ \displaystyle{\sup_{j\in \mathbb{Z}}\| \varphi_{j}\|_{1,k}<+\infty}$;
\item[$(iii)$] $\displaystyle{\sum_{j\in \mathbb{Z}}\mathcal{F}_k(\varphi_{j})(x) =1}$,  for  $x\in \mathbb{R}^d\backslash\{0\}$.
\end{itemize}
$ \mathcal{S}(\mathbb{R}^d)^{\rm rad}$ being the subspace of functions
in the Schwartz space $\mathcal{S}(\mathbb{R}^d)$ which are
radial.

Put $\displaystyle{\mathcal{A}=\big\{ \phi \in
\mathcal{S}(\mathbb{R}^d)^{\rm rad}: \,\mbox{supp}\,\mathcal{F}_k(\phi)
\subset \{x\in \mathbb{R}^d ;\,1\leq \|x\|\leq 2\} \big\}}.$ Given
$\phi \in \mathcal{A}$, we denote by~$ \mathcal{C}_{ p,q}^{\phi
,\beta,k} $ the subspace of functions $f \in L^p_k(\mathbb{R}^d)$
satisfying
\begin{gather*}
\bigg(\int^{+ \infty}_0 \left(\frac{\|f \ast_k \phi_t
\|_{p,k}}{t^\beta}\right)^q \, \frac{dt}{t}\bigg)^{\frac{1}{q}}  <  +
\infty \qquad \mbox{if}\quad q < +\infty
\end{gather*} and
\begin{gather*}  \sup_{t\in (0,+\infty)}\frac{\|f \ast_k \phi_t\|_{p,k}}{t^\beta}  < +\infty \qquad \mbox{if} \quad q = +\infty ,\end{gather*}
where
$\phi_t(x)=\frac{1}{t^{2(\gamma+\frac{d}{2})}}\phi(\frac{x}{t})$,
 for all $t\in (0,+\infty)$ and $x\in\mathbb{R}^d$.

In this paper we show for $ \beta>0$ the inclusion of the Schwartz
space in $\mathcal{B}\mathcal{D}_{p,q}^{\beta,k}$ for $1\leq
p,q\leq+\infty$ and the density when $1\leq p,q<+\infty$. We prove
an interpolation formula for the Besov--Dunkl spaces by the real
method. We compare these spaces with $\mathcal{C}_{ p,q}^{\phi
,\beta,k}$ which extend to the Dunkl operators on $\mathbb{R}^d$
some results obtained in \cite{[4], [5], [21]}. Finally we establish further
results of integrability of $\mathcal{F}_k(f) $ when $f$ is in a
suitable
 Besov--Dunkl space $\mathcal{B}\mathcal{D}_{p,q}^{\beta,k} $ for $1\leq p\leq2$ and $1\leq q\leq+\infty$.
 Using the characterization of the Besov spaces by dif\/ferences analogous results of integrability have been obtained in the case $q=1$ by Giang and
M\'{o}ricz in~\cite{[13]} for a~classical Fourier transform on~$\mathbb{R}$
 and for $q=1,\,+\infty$ by Betancor and Rodr\'{\i}guez-Mesa in~\cite{[7]} for the Hankel transform on $(0,+\infty)$ in Lipschitz--Hankel spaces. Later
 Abdelkef\/i and Sif\/i in~\cite{[1], [2]} have established similar results of
 integrability for the Dunkl transform on $\mathbb{R}$ and in radial case
 on~$\mathbb{R}^d$. The argument used in~\cite{[1], [2]} to establish such integrability
 is the $L^p$-boundedness of the Dunkl translation operators, making it
 dif\/f\/icult to extend the results on $\mathbb{R}^d$. We take a
 dif\/ferent approach based on the the characterization of the Besov spaces by convolution to establish
 our results on higher dimension.

The contents of this paper are as follows. In Section~\ref{section2} we
collect some basic def\/initions and results about harmonic analysis
associated with Dunkl operators.
In Section~\ref{section3} we show the inclusion and the density of the Schwartz
space in $\mathcal{B}\mathcal{D}_{p,q}^{\beta,k}$, we prove an
interpolation formula for the Besov--Dunkl spaces by the real method
and we compare these spaces with $\mathcal{C}_{ p,q}^{\phi
,\beta,k}$.
 In Section~\ref{section4} we establish our results concerning
integrability of the Dunkl transform of function in the Besov--Dunkl
spaces.

Along this paper we denote by $\langle \cdot,\cdot\rangle$ the usual
Euclidean inner product in $\mathbb{R}^d$ as well as its extension
to $\mathbb{C}^d \times \mathbb{C}^d$, we write for $x \in
\mathbb{R}^d, \|x\| = \sqrt{\langle x,x\rangle}$ and we represent by
$c$ a suitable positive constant which is not necessarily the same
in each occurrence. Furthermore we denote~by
\begin{itemize}\itemsep=0pt
\item $\mathcal{E}(\mathbb{R}^d)$ the space of inf\/initely
dif\/ferentiable functions on $\mathbb{R}^d$;

\item $\mathcal{S}(\mathbb{R}^d)$ the Schwartz space of
functions in $\mathcal{E}( \mathbb{R}^d)$ which are rapidly
decreasing as well as their derivatives;

\item $\mathcal{D}(\mathbb{R}^d)$ the subspace of
$\mathcal{E}(\mathbb{R}^d)$ of compactly supported functions.
\end{itemize}

\section{Preliminaries}\label{section2}

Let $W$ be a f\/inite ref\/lection group on $\mathbb{R}^{d}$, associated
with a root system $R$ and $R_+$ the positive subsystem of $R$ (see
\cite{[8], [10], [11], [12], [18], [19]}). We denote by $k$ a nonnegative
multiplicity function def\/ined on $R$ with the property that $k$ is
$W$-invariant. We associate with $k$ the index
\[
\gamma = \gamma (R) = \sum_{\xi \in R_+} k(\xi) \geq 0,
\]
and the weight function $w_k$ def\/ined by
\[
w_k(x) = \prod_{\xi\in R_+} |\langle \xi,x\rangle|^{2k(\xi)},
\qquad x \in \mathbb{R}^d .
\]  Further we introduce the Mehta-type
constant $c_k$ by
\[
c_k = \left(\int_{\mathbb{R}^d} e^{- \frac{\|x\|^2}{2}}
w_k (x)dx\right)^{-1}.
\]

For every $1 \leq p \leq + \infty$ we denote by
$L^p_k(\mathbb{R}^d)$ the space $L^{p}(\mathbb{R}^d, w_k(x)dx),$
 $L^p_k( \mathbb{R}^d)^{\rm rad}$  the subspace of those $f \in
L^p_k( \mathbb{R}^d)$ that are radial  and  we use $\|\cdot \|_{p,k}$
as a shorthand for $\|\cdot \|_{L^p_k( \mathbb{R}^d)}.$

By using the homogeneity of $w_k$ it is shown in~\cite{[18]} that for $f
\in L^1_k ( \mathbb{R}^d)^{\rm rad}$ there exists a~function $F$ on
$[0, + \infty)$ such that $f(x) = F(\|x\|)$, for all $x \in
\mathbb{R}^d$. The function $F$ is integrable with respect to the
measure $r^{2\gamma+d-1}dr$ on $[0, + \infty)$ and we have
 \begin{gather*}\int_{S^{d-1}}w_k (x)d\sigma(x) = \frac{c^{-1}_k}{2^{\gamma +\frac{d}{2} -1}
\Gamma(\gamma + \frac{d}{2})} ,   \end{gather*} where $S^{d-1}$
is the unit sphere on $\mathbb{R}^d$ with the normalized surface
measure $d\sigma$  and
 \begin{gather} \int_{\mathbb{R}^d}  f(x)w_k(x)dx  =
 \int^{+\infty}_0
\Big( \int_{S^{d-1}}w_k(ry)d\sigma(y)\Big)
F(r)r^{d-1}dr\nonumber\\
\phantom{\int_{\mathbb{R}^d}  f(x)w_k(x)dx}{}  =   \frac{c^{-1}_k}{2^{\gamma +
\frac{d}{2}-1}\Gamma(\gamma + \frac{d}{2})}\int^{+ \infty}_0 F(r)
r^{2\gamma+d-1}dr.\label{eq1}
\end{gather}

Introduced by C.F.~Dunkl in~\cite{[9]} the Dunkl operators $T_j$, $1\leq j\leq d$, on $\mathbb{R}^d$ associated with the ref\/lection
group $W$ and the multiplicity function $k$ are the f\/irst-order
dif\/ferential- dif\/ference operators given by
\[
T_jf(x)=\frac{\partial f}{\partial x_j}(x)+\sum_{\alpha\in\mathcal{R}_+}k(\alpha)
\alpha_j\,\frac{f(x)-f(\sigma_\alpha(x))}{\langle\alpha,x\rangle} ,\qquad
f\in\mathcal{E}(\mathbb{R}^d),\qquad x\in\mathbb{R}^d,
\] where
$\alpha_j=\langle\alpha,e_j\rangle$, $(e_1,e_2,\ldots,e_d)$ being
the canonical basis of $\mathbb{R}^d$.

The Dunkl kernel $E_k $ on $\mathbb{R}^d \times \mathbb{R}^d$ has
been introduced by C.F.~Dunkl in~\cite{[10]}. For $y \in \mathbb{R}^d$ the
function $x \mapsto E_k(x,y)$ can be viewed as the solution on
$\mathbb{R}^d$ of the following initial problem
\begin{gather*}
T_ju(x,y) = y_j\,u(x,y),\qquad 1\leq j\leq d,\\
u(0,y)=1.
\end{gather*}  This kernel has a unique
holomorphic extension to $\mathbb{C}^d \times \mathbb{C}^d $. M.~R\"osler has proved in~\cite{[17]} the following integral representation
for the Dunkl kernel
\[
E_k(x,z)=\int_{\mathbb{R}^d}e^{\langle
y,z\rangle}d\mu_x^k(y),\qquad x\in\mathbb{R}^d,\qquad
z\in\mathbb{C}^d,
\] where $\mu_x^k$ is a probability measure on $\mathbb{R}^d$ with support in the closed ball
$B(0,\|x\|)$ of center $0$ and radius $\|x\|$.

We have for all $\lambda\in \mathbb{C}$ and $z, z'\in
\mathbb{C}^d$
$E_k(z,z') = E_k(z',z)$,  $E_k(\lambda z,z') = E_k(z,\lambda z')$ and for $x, y
\in \mathbb{R}^d$
 $|E_k(x,iy)| \leq 1$ (see \cite{[10], [17], [18], [19], [22]}).

The Dunkl transform $\mathcal{F}_k$ which was introduced by C.F.~Dunkl
in~\cite{[11]} (see also~\cite{[8]}) is def\/ined for $f \in \mathcal{D}(
\mathbb{R}^d)$ by
\[
\mathcal{F}_k(f)(x) =c_k\int_{\mathbb{R}^d}f(y) E_k(-ix, y)w_k(y)dy,\qquad
x \in \mathbb{R}^d.
\] According to \cite{[8], [11], [18]} we have the
following results:
\begin{itemize}\itemsep=0pt
\item[$i)$] The Dunkl transform of a function $f
\in L^1_k( \mathbb{R}^d)$ has the following basic property
\begin{gather}\label{eq2}
\| \mathcal{F}_k(f)\|_{\infty,k} \leq
 \|f\|_{ 1,k}.
  \end{gather}
\item[$ii)$] The Schwartz space $\mathcal{S}( \mathbb{R}^d)$ is
invariant under the Dunkl transform $\mathcal{F}_k\;.$
\item[$iii)$] When both $f$ and $\mathcal{F}_k(f)$ are in $L^1_k( \mathbb{R}^d)$,
 we have the inversion formula
 \begin{gather*} f(x) =   \int_{\mathbb{R}^d}\mathcal{F}_k(f)(y) E_k( ix, y)w_k(y)dy,\qquad
x \in \mathbb{R}^d.
\end{gather*}
\item[$iv)$] (Plancherel's theorem) The Dunkl transform on $\mathcal{S}(\mathbb{R}^d)$
 extends uniquely to an isometric isomorphism on
$L^2_k(\mathbb{R}^d)$.
\end{itemize}

By \eqref{eq2}, Plancherel's theorem and the
Marcinkiewicz interpolation theorem (see~\cite{[20]}) we get for $f \in
L^p_k(\mathbb{R}^d)$ with $1\leq p\leq 2$ and $p'$ such that
$\frac{1}{p}+\frac{1}{p'}=1$,
\begin{gather}\label{eq3}
\|\mathcal{F}_k(f)\|_{p',k}\leq c \|f\|_{p,k}.
\end{gather}
The Dunkl transform of a function in $L^1_k( \mathbb{R}^d)^{\rm rad}$
 is also radial and could be expressed via the Hankel transform (see~\cite[Proposition~2.4]{[18]}).

K.~Trim\`eche has introduced in~\cite{[23]} the Dunkl translation operators
$\tau_x$, $x\in\mathbb{R}^d$, on $\mathcal{E}( \mathbb{R}^d)$. For
$f\in \mathcal{S}( \mathbb{R}^d)$ and $x, y\in\mathbb{R}^d$ we have
\[
\mathcal{F}_k(\tau_x(f))(y)=E_k(i x, y)\mathcal{F}_k(f)(y)
\] and
\begin{gather}\label{eq4}
\tau_x(f)(y)=c_k \int_{\mathbb{R}^d}\mathcal{F}_k (f)(\xi)
E_k(ix,\xi)E_k(iy,\xi)w_k(\xi)d\xi.
\end{gather} Notice that for all $x,y\in\mathbb{R}^d$
$\tau_x(f)(y)=\tau_y(f)(x)$, and for f\/ixed $x\in\mathbb{R}^d$
\begin{gather}\label{eq5}
\tau_x \ \mbox{is a continuous linear mapping from} \
\mathcal{E}( \mathbb{R}^d) \ \mbox{into} \ \mathcal{E}(\mathbb{R}^d).
\end{gather} As an operator on
$L_k^2(\mathbb{R}^d)$, $\tau_x$ is bounded. A priori it is not at
all clear whether the translation operator can be def\/ined for $L^p$-functions with $p$ dif\/ferent from~2. However, according to~\cite[Theorem~3.7]{[19]} the operator $\tau_x$ can be extended to
$L^p_k(\mathbb{R}^d)^{\rm rad}$, $1 \leq p \leq 2$ and for $f \in
L^p_k(\mathbb{R}^d)^{\rm rad}$ we have
\[
\|\tau_x(f)\|_{p,k} \leq \|f\|_{p,k}.
\]

The Dunkl convolution product $\ast_k$ of two functions $f$ and $g$
in $L^2_k(\mathbb{R}^d)$ (see~\cite{[19], [23]}) is given~by
\[
(f\; \ast_k g)(x) = \int_{\mathbb{R}^d} \tau_x (f)(-y) g(y) w_k(y)dy,\qquad
x \in \mathbb{R}^d .
\] The Dunkl convolution product is commutative
and for $f, g \in \mathcal{D}( \mathbb{R}^d)$ we have
\begin{gather}\label{eq6}
\mathcal{F}_k(f\,\ast_k\, g) =
\mathcal{F}_k(f) \mathcal{F}_k(g).
\end{gather}
It was shown in~\cite[Theorem 4.1]{[19]} that when $g$ is a bounded function in $L^1_k(
\mathbb{R}^d)^{\rm rad}$, then
\[
(f  \ast_k g)(x) = \int_{\mathbb{R}^d}  f(y) \tau_x (g)(-y) w_k(y)dy,\qquad
x \in \mathbb{R}^d ,
\] initially def\/ined on the intersection of
$L^1_k(\mathbb{R}^d)$ and $L^2_k(\mathbb{R}^d)$ extends to all
$L^p_k(\mathbb{R}^d)$, $1\leq p\leq +\infty$ as a bounded operator.
In particular,
\begin{gather}\label{eq7}
\|f \ast_k g\|_{p,k} \leq \|f\|_{p,k}
\|g\|_{1,k}.
\end{gather}

The Dunkl Laplacian $\Delta_k$ is def\/ined by
$\Delta_k:=\sum\limits_{i=1}^{d}T_i^2 $. From~\cite{[16]} we
have for each $\lambda>0$ $\lambda I-\Delta_k$ maps
$\mathcal{S}(\mathbb{R}^d)$ onto itself  and
\begin{gather}\label{eq8}
\mathcal{F}_k((\lambda
I-\Delta_k)f)(x)=\big(\lambda+\|x\|^2\big)\mathcal{F}_k(f)(x), \qquad \mbox{for}\quad x\in\mathbb{R}^d.
\end{gather}

\section[Interpolation and characterization for the Besov-Dunkl spaces]{Interpolation and characterization\\ for the Besov--Dunkl spaces}\label{section3}

In this section we establish the inclusion and the density of
$\mathcal{S}(\mathbb{R}^d)$ in
$\mathcal{B}\mathcal{D}_{p,q}^{\beta,k}$ and we prove an
interpolation formula for the Besov--Dunkl spaces by the real method.
Finally we compare the spaces
$\mathcal{B}\mathcal{D}_{p,q}^{\beta,k}$ with $\mathcal{C}_{
p,q}^{\phi ,\beta,k}$. Before, we start with some useful
results.

We shall denote by $\mathbf{\Phi}$ the set of all sequences of
functions $ (\varphi_{j})_{j\in \mathbb{Z}}$ in $
\mathcal{S}(\mathbb{R}^d)^{\rm rad}$ satisfying
\begin{itemize}\itemsep=0pt
\item[$(i)$] $\mbox{supp}\,\mathcal{F}_k(\varphi_{j}) \subset A_j=\{x\in
\mathbb{R}^d ;\,2^{j-1}\leq \|x\|\leq 2^{j+1}\}$ for $j\in
\mathbb{Z}$;
\item[$(ii)$] $ {\sup\limits_{j\in \mathbb{Z}}\| \varphi_{j}\|_{1,k}<+\infty}$;
\item[$(iii$)] ${\sum\limits_{j\in \mathbb{Z}}\mathcal{F}_k(\varphi_{j})(x) =1}$,  for  $x\in \mathbb{R}^d\backslash\{0\}$.
\end{itemize}

\begin{proposition}\label{proposition3.1} Let $ \beta>0$ and $1\leq p , q \leq +
\infty$, then $\mathcal{B}\mathcal{D}_{p,q}^{\beta,k}$ is
independent of the choice of the sequence in~$\mathbf{\Phi}$.
\end{proposition}

\begin{proof} Fix $ (\phi_{j})_{j\in \mathbb{Z}}$, $ (
\psi_{j})_{j\in \mathbb{Z}}$ in $\mathbf{\Phi}$ and
$f\in\mathcal{B}\mathcal{D}_{p,q}^{\beta,k}$ for $q<+\infty$. Using
the properties
 $(i)$ for $ (\phi_{j})_{j\in \mathbb{Z}}$ and $(i)$
and $(iii)$ for $ ( \psi_{j})_{j\in \mathbb{Z}}$, we have for $j\in
\mathbb{Z}$ $\phi_j= \phi_j\ast_k(
\psi_{j-1}+\psi_{j}+\psi_{j+1})$. Then by the property $(ii)$ for $
(\phi_{j})_{j\in \mathbb{Z}}$, \eqref{eq7} and H\"older's inequality for
$j\in \mathbb{Z}$ we obtain
\begin{gather*} \|f \ast_k \phi_j \|_{p,k}^q\leq c\, 3^{q-1}
\sum_{s=j-1}^{j+1}\| \psi_s\ast_k f\|_{p,k}^q.
\end{gather*} Thus summing over $j$ with weights $2^{j \beta q}$
we get \begin{gather*}
 \sum_{j\in\mathbb{Z}}\big(2^{j\beta}\|\phi_{j }\ast_k f\|_{p,k}\big)^q\leq
 c \sum_{j\in\mathbb{Z}}\big(2^{j\beta}\|\psi_{j }\ast_k f\|_{p,k}\big)^q.
   \end{gather*} Hence by symmetry we get the result of our proposition.
   When $q=+\infty$ we make the usual modif\/ication.
 \end{proof}

 \begin{remark}\label{remark3.1}
Let $\beta>0$,
$1\leq p , q \leq + \infty$, we denote by $
\mathcal{\ddot{B}}\mathcal{D}_{p,q}^{\beta,k} $ the subspace of
functions $f \in L^p_k(\mathbb{R}^d)$ satisfying
\[
\bigg(\sum_{j\in\mathbb{N}}(2^{j\beta}\|\varphi_{j }\ast_k
f\|_{p,k})^q \bigg)^{\frac{1}{q}}<+\infty \qquad \mbox{if}\quad q <
+\infty
\] and
\[
 \sup_{j\in \mathbb{N}}2^{j\beta}\|\varphi_{j
}\ast_k f\|_{p,k}<+\infty\qquad \mbox{if} \quad q = +\infty ,
\] where $
(\varphi_{j})_{j\in \mathbb{N}}$ is a sequence of functions in $
\mathcal{S}(\mathbb{R}^d)^{\rm rad}$ such that
\begin{itemize}\itemsep=0pt
\item[$i)$] $\mbox{supp}\,\mathcal{F}_k(\varphi_{0})
\subset \{x\in \mathbb{R}^d ;\, \|x\|\leq 2\}$  and
$\mbox{supp}\,\mathcal{F}_k(\varphi_{j}) \subset A_j=\{x\in
\mathbb{R}^d ;\,2^{j-1}\leq \|x\|\leq 2^{j+1}\}$ for $j\in
\mathbb{N\backslash}\{0\}$;
\item[$ii)$] $ {\sup\limits_{j\in \mathbb{N}}\| \varphi_{j}\|_{1,k}<+\infty}$;
\item[$iii$)] $ {\sum\limits_{j\in \mathbb{N}}\mathcal{F}_k(\varphi_{j})(x) =1}$,  for  $x\in \mathbb{R}^d\backslash\{0\}$.
\end{itemize}
As the Besov--Dunkl spaces these spaces are also independent of the choice of the sequence $
(\varphi_{j})_{j\in \mathbb{N}}$ satisfying the previous properties.
\end{remark}

\begin{proposition}\label{proposition3.2} For $ \beta>0$ and $1\leq p , q \leq + \infty$ we have
\[
\mathcal{\ddot{B}}\mathcal{D}_{p,q}^{\beta,k}=\mathcal{
B}\mathcal{D}_{p,q}^{\beta,k}.
\]
\end{proposition}

\begin{proof} Since both spaces $
\mathcal{\ddot{B}}\mathcal{D}_{p,q}^{\beta,k}$ and $\mathcal{
B}\mathcal{D}_{p,q}^{\beta,k}$ are in $ L^p_k(\mathbb{R}^d)$ and are
independent of the specif\/ic selection of sequence of functions, then
according to \cite[Lemma~6.1.7, Theorem~6.3.2]{[5]} we can take a function $\phi \in
\mathcal{S}(\mathbb{R}^d)^{\rm rad}$ such that
\begin{itemize}\itemsep=0pt
\item $\mbox{supp}\,\mathcal{F}_k(\phi)\subset \{x\in \mathbb{R}^d ;\, \frac{1}{2}\leq \|x\|\leq 2\}$;
\item
$\mathcal{F}_k(\phi)(x) > 0$ for $\frac{1}{2}< \|x\|<2$;

\item $ {\sum\limits_{j\in \mathbb{Z}}\mathcal{F}_k(\phi_{2^{-j}})(x) =1}$, $x\in \mathbb{R}^d\backslash\{0\}$.
\end{itemize}
If we consider the sequences $(\psi_{j})_{j\in \mathbb{Z}}$ and
$(\varphi_{j})_{j\in \mathbb{N}}$ in
$\mathcal{S}(\mathbb{R}^d)^{\rm rad}$ def\/ined respectively for
$\mathcal{ B}\mathcal{D}_{p,q}^{\beta,k}$ and
$\mathcal{\ddot{B}}\mathcal{D}_{p,q}^{\beta,k}$ by
$\psi_{j}=\phi_{2^{-j}}$
 $\forall \,j\in \mathbb{Z}$ and $\varphi_{0}= {\sum\limits_{j\in
\mathbb{Z}_-}\phi_{2^{-j}}}$, $\varphi_{j}= \phi_{2^{-j}}$ $
\forall \,j\in \mathbb{N}^*$, we can assert that $
\mathcal{\ddot{B}}\mathcal{D}_{p,q}^{\beta,k}=\mathcal{
B}\mathcal{D}_{p,q}^{\beta,k}$.
\end{proof}

\begin{remark}\label{remark3.2}
By Proposition~\ref{proposition3.2} and~\cite[Proposition~2]{[21]} we have the following embeddings.
\begin{enumerate}\itemsep=0pt
\item[1.] Let $1\leq q_1 \leq q_2 \leq +\infty$ and $\beta >0$.
Then
\[
\mathcal{B}\mathcal{D}_{p,q_1}^{\beta,k} \subset
\mathcal{B}\mathcal{D}_{p,q_2}^{\beta,k}\qquad \mbox{if}\quad 1\leq
p \leq +\infty.
\]
\item[2.] Let $1\leq q_1 , q_2 \leq +\infty$, $\beta >0$ and $\varepsilon
>0$. Then
\[
\mathcal{B}\mathcal{D}_{p,q_1}^{\beta+\varepsilon ,k} \subset
\mathcal{B}\mathcal{D}_{p,q_2}^{\beta,k}\qquad \mbox{if}\quad 1\leq
p \leq +\infty.
\]
\end{enumerate}
\end{remark}

\begin{proposition}\label{proposition3.3} For $\beta>0$ and $1\leq p , q \leq +
\infty$ we have
\[
\mathcal{S}(\mathbb{R}^d)\subset
\mathcal{B}\mathcal{D}_{p,q}^{\beta,k}.
\]
 If $1\leq p , q < +
\infty$, then $\mathcal{S}(\mathbb{R}^d)$ is dense in
$\mathcal{B}\mathcal{D}_{p,q}^{\beta,k}.$
\end{proposition}

\begin{proof} In order to prove the inclusion, we may restrict
ourself to $q=+\infty$. This follows from the fact that
$\mathcal{B}\mathcal{D}_{p,\infty}^{\beta ,k} \subset
\mathcal{B}\mathcal{D}_{p,q}^{\beta',k}$ for $\beta>\beta'>0$ and
$1\leq p,q \leq +\infty$ (see Remark~\ref{remark3.2},~2). Let $f\in
\mathcal{S}(\mathbb{R}^d)$ and $ (\varphi_{j})_{j\in \mathbb{N}}$ a
sequence of functions in $\mathcal{S}(\mathbb{R}^d)^{\rm rad}$
satisfying the properties of Remark~\ref{remark3.1}, there exists a suf\/f\/iciently
large natural number $L$ such that
\[
 \sup_{j\in\mathbb{N}}2^{j\beta}\|\varphi_j\ast_k f\|_{p,k} \leq  \sup_{j\in
\mathbb{N}}2^{j\beta}\|(1+\|x\|^2)^L(\varphi_j\ast_k
f)\|_{\infty,k}.
\] Since
$\varphi_j\in\mathcal{S}(\mathbb{R}^d)^{\rm rad}$, then for $x\in
\mathbb{R}^d\,$
$\mathcal{F}_k^{-1}(\varphi_j)(x)=\mathcal{F}_k(\varphi_j)(-x)=\mathcal{F}_k(\varphi_j)(x)$,
so using \eqref{eq8} and the property i) of $ (\varphi_{j})_{j\in
\mathbb{N}}$ (see Remark~\ref{remark3.1}) we
obtain\begin{gather*}
\sup_{j\in\mathbb{N}^{*}}2^{j\beta}\|\varphi_j\ast_k
f \|_{p,k} \leq \sup_{j\in
\mathbb{N}^{*}}2^{j\beta}\|\mathcal{F}_k[(I-\Delta_k)^L(\mathcal{F}_k(\varphi_j)\mathcal{F}_k^{-1}(f))] \|_{\infty,k}\\
\phantom{\sup_{j\in\mathbb{N}^{*}}2^{j\beta}\|\varphi_j\ast_k
f \|_{p,k}}{} \leq \sup_{j\in
\mathbb{N}^{*}}2^{j\beta}\|(I-\Delta_k)^L(\mathcal{F}_k(\varphi_j)\mathcal{F}_k^{-1}(f))\|_{1,k}\\
\phantom{\sup_{j\in\mathbb{N}^{*}}2^{j\beta}\|\varphi_j\ast_k
f \|_{p,k}}{} \leq \sup_{j\in \mathbb{N}^{*}}\,c_j \,2^{j\beta}\sup_{x\in
A_j}|(I-\Delta_k)^L(\mathcal{F}_k(\varphi_j)\mathcal{F}_k^{-1}(f))(x)|,
\end{gather*}
where
 $c_j= {\int_{A_j}w_k(x)} dx$. Hence there exists a suf\/f\/iciently large
natural number $M$ such that $ {\frac{c_j
2^{j\beta}}{(1+2^{2(j-1)})^M}}\leq 1$, $\forall\, j\in \mathbb{N}^{*}$
and we get
\begin{gather*}
\sup_{j\in\mathbb{N}^{*}}2^{j\beta}\|\varphi_j\ast_k f\|_{p,k} \leq
 \sup_{j\in \mathbb{N}^{*}}\sup_{x\in
A_j}|\big(1+\|x\|^2\big)^M(I-\Delta_k)^L(\mathcal{F}_k(\varphi_j)\mathcal{F}_k^{-1}(f))(x)|.
\end{gather*}
Since $(I-\Delta_k)^L$ is linear and continuous from
$\mathcal{S}(\mathbb{R}^d)$ into itself, we deduce that
\begin{gather*}\sup_{j\in\mathbb{N}^{*}}2^{j\beta}\|\varphi_j\ast_k f\|_{p,k}\leq
 c \sup_{j\in \mathbb{N}^{*}}\sup_{x\in
\mathbb{R}^d}|\mathcal{F}_k(\varphi_j)(x)|\sup_{x\in
\mathbb{R}^d}|(1+\|x\|^2)^M\mathcal{F}_k^{-1}(f)(x)|,
\end{gather*}
which gives by the property ii) of $ (\varphi_{j})_{j\in
\mathbb{N}}$
\begin{gather*}\sup_{j\in\mathbb{N}^{*}}2^{j\beta}\|\varphi_j\ast_k
f\|_{p,k}\leq c \sup_{j\in \mathbb{N}^{*}}\|
\varphi_j\|_{1,k}\sup_{x\in
\mathbb{R}^d}|(1+\|x\|^2)^M\mathcal{F}_k^{-1}(f)(x)|<+\infty.
\end{gather*}
By Proposition~\ref{proposition3.2} we conclude that
$\mathcal{S}(\mathbb{R}^d)\subset
\mathcal{B}\mathcal{D}_{p,q}^{\beta,k}.$

Let us now prove the
density of $\mathcal{S}(\mathbb{R}^d)$ in
$\mathcal{B}\mathcal{D}_{p,q}^{\beta,k}$ for $p,q<+\infty$. Assume
$f\in\mathcal{B}\mathcal{D}_{p,q}^{\beta,k}$ and $
(\varphi_{j})_{j\in \mathbb{Z}}\in\mathbf{\Phi}$, then we put for
$N\in\mathbb{N\backslash}\{0\}$, $f_N= { \sum\limits_{s=-N}^{N}
\varphi_s\ast_k f}.$ It's clear that
$f_N\in\mathcal{B}\mathcal{D}_{p,q}^{\beta,k}$. We have
\begin{gather*}
\sum_{j\in\mathbb{Z}}2^{j\beta q}\|\varphi_{j }\ast_k(
f_N-f)\|_{p,k}^q=\sum_{j\in\mathbb{Z}}2^{j\beta q}\left\|\bigg(\sum_{s=-N}^{N}
\varphi_s\bigg)\ast_k\varphi_{j }\ast_k f- \varphi_j\ast_k
f \right\|_{p,k}^q.
\end{gather*} Using the properties $(i)$ and $(iii)$ for $
( \varphi_{j})_{j\in \mathbb{Z}}$ we get
\[
\|f-f_N\|^q_{\mathcal{B}\mathcal{D}_{p,q}^{\beta,k}}\leq c \sum_{|j|\geq N}
2^{j\beta q}\|\varphi_j\ast_k f\|_{p,k}^q.
\]
 Since
$f\in\mathcal{B}\mathcal{D}_{p,q}^{\beta,k}$, then we deduce that
\begin{gather}\label{eq9}
\lim_{N\rightarrow+\infty}\|f-f_N\|_{\mathcal{B}\mathcal{D}_{p,q}^{\beta,k}}=0.
\end{gather}
Next we take a function $\theta\in \mathcal{D}(\mathbb{R}^d)$ such
that $\theta(0)=1$. For $n\in\mathbb{N}\backslash\{0\}$ we put
$\theta_n(x)=\theta(n^{-1}x)$, $x\in\mathbb{R}^d$. From~\eqref{eq5} we have
for $N\in\mathbb{N}\backslash\{0\}$
$f_N\in\mathcal{E}(\mathbb{R}^d)$, then
$f_N\theta_n\in\mathcal{S}(\mathbb{R}^d)$. Again from the properties
$(i)$ and $(iii)$ for $ ( \varphi_{j})_{j\in \mathbb{Z}}$ we can assert
that $f_N= { \sum\limits_{j=-N}^{N} \varphi_j\ast_k f_{N+1}}$
which gives $f_N\theta_n= { \sum\limits_{j=-N}^{N}
\varphi_j\ast_k f_{N+1}\theta_n}$. Using the properties $(i)$, $(ii)$,
$(iii)$ for $ ( \varphi_{j})_{j\in \mathbb{Z}}$ and \eqref{eq7} we obtain
\begin{gather*}
\|f_N-f_N\theta_n\|^q_{\mathcal{B}\mathcal{D}_{p,q}^{\beta,k}}\leq c
\sum_{j=-N-1}^{N+1} 2^{j\beta q}\|f_{N+1}-f_{N+1}\theta_n\|_{p,k}^q.
\end{gather*} The dominated convergence theorem implies that
\[
\|f_{N+1}-f_{N+1}\theta_n\|_{p,k} \rightarrow 0\qquad \mbox{as}\quad
 n\rightarrow +\infty .
 \]  Hence we deduce that
  \begin{gather}\label{eq10}
  \|f_N-f_N\theta_n\|_{\mathcal{B}\mathcal{D}_{p,q}^{\beta,k}} \rightarrow 0 \qquad \mbox{as} \quad n \rightarrow
+\infty, \end{gather}
Combining \eqref{eq9} and \eqref{eq10} we conclude that
$\mathcal{S}(\mathbb{R}^d)$ is dense in
$\mathcal{B}\mathcal{D}_{p,q}^{\beta,k}$. This completes the proof
of Proposition~\ref{proposition3.3}.
\end{proof}

For $0<\theta<1$, $\beta_0, \beta_1
>0$, $1\leq p, q_0,q_1\leq +\infty$ and $1\leq q\leq +\infty$, the
real interpolation Besov--Dunkl space denoted by
$(\mathcal{B}\mathcal{D}_{p,q_0}^{\beta_0,k},
\mathcal{B}\mathcal{D}_{p,q_1}^{\beta_1,k})_{\theta,q}$ is the
subspace of functions $f\in
\mathcal{B}\mathcal{D}_{p,q_0}^{\beta_0,k} +
 \mathcal{B}\mathcal{D}_{p,q_1}^{\beta_1,k}$ satisfying
\[
\left(\int^{+ \infty}_0 \left(t^{-\theta}\mathcal{K}_{p,k}(t,f; \beta_0,q_0 ; \beta_1,q_1)\right)^q   \frac{dt}{t}\right)^{\frac{1}{q}} < + \infty
\qquad\mbox{if} \quad q <+\infty ,
\]
 and
 \[
 \sup_{t\in
(0,+\infty)}t^{-\theta}\mathcal{K}_{p,k}(t,f; \beta_0,q_0 ; \beta_1,q_1)
<+\infty \qquad\mbox{if} \quad q = +\infty ,
\] with
$\mathcal{K}_{p,k}$ is the Peetre $\mathcal{K}$-functional given by
\[
\mathcal{K}_{p,k}(t,f; \beta_0,q_0; \beta_1,q_1)=\inf\Big\{\| f_0\|_{\mathcal{B}\mathcal{D}_{p,q_0}^{\beta_0,k}}+t\|f_1\|_{\mathcal{B}\mathcal{D}_{p,q_1}^{\beta_1,k}}\Big\},
\]
 where the inf\/inimum is taken over all representations of $f$ of the form
 \[
 f=f_0+f_1, \qquad f_0\in\mathcal{B}\mathcal{D}_{p,q_0}^{\beta_0,k},\qquad
  f_1\in\mathcal{B}\mathcal{D}_{p,q_1}^{\beta_1,k}.
 \]

\begin{theorem}\label{theorem3.1}
Let $0<\theta<1$ and $1\leq p,q,q_0,q_1\leq
+\infty$. For $\beta_0, \beta_1 >0$, $\beta_0\neq \beta_1$ and
$\beta=(1-\theta)\beta_0+ \theta\beta_1$ we have
\[
\big(\mathcal{B}\mathcal{D}_{p,q_0}^{\beta_0,k}, \mathcal{B}\mathcal{D}_{p,q_1}^{\beta_1,k}\big)_{\theta,q}=
\mathcal{B}\mathcal{D}_{p,q}^{\beta,k}.
\]
\end{theorem}

\begin{proof} We start with
the proof of the inclusion
$(\mathcal{B}\mathcal{D}_{p,\infty}^{\beta_0,k}, \mathcal{B}
\mathcal{D}_{p,\infty}^{\beta_1,k})_{\theta,q}\subset
\mathcal{B}\mathcal{D}_{p,q}^{\beta,k}.$
We may assume that $\beta_0>\beta_1$. Let $q<+\infty$, for
$f=f_0+f_1$ with $f_0\in \mathcal{B}\mathcal{D}_
{p,\infty}^{\beta_0,k}$ and
$f_1\in\mathcal{B}\mathcal{D}_{p,\infty}^{\beta_1,k}$ we get by
Proposition~\ref{proposition3.2}
\begin{gather*}
\sum _{l=0}^{+\infty}2^{ql\beta}\|\varphi_l \ast_k f\|_{p,k}^q \leq
c\sum _{l=0}^{+\infty}2^{-\theta ql(\beta_0-\beta_1)}
\Big(2^{l\beta_0}\|\varphi_l \ast_k f_0\|_{p,k}+
2^{l(\beta_0-\beta_1)}2^{l\beta_1}\|\varphi_l \ast_k
f_1\|_{p,k}\Big)^q  \\
\phantom{\sum _{l=0}^{+\infty}2^{ql\beta}\|\varphi_l \ast_k f\|_{p,k}^q}{}  \leq   c\sum _{l=0}^{+\infty}2^{-\theta q
l(\beta_0-\beta_1)}\Big(\|f_0\|_{\mathcal{B}\mathcal{D}_{p,\infty}^{\beta_0,k}}
+ 2^{l(\beta_0-\beta_1)}\|f_1 \|_{\mathcal{
B}\mathcal{D}_{p,\infty}^{\beta_1,k}}\Big)^q.
\end{gather*}
Then we deduce that
\begin{gather*}\sum _{l=0}^{+\infty}2^{ql\beta}\|\varphi_l
\ast_k f\|_{p,k}^q  \leq  c  \sum _{l=0}^{+\infty}2^{-\theta q
l(\beta_0-\beta_1)}\Big(\mathcal{K}_{p,k}(2^{l(\beta_0-\beta_1)}
,f; \beta_0,\infty ; \beta_1,\infty)\Big)^q
 \\
 \phantom{\sum _{l=0}^{+\infty}2^{ql\beta}\|\varphi_l
\ast_k f\|_{p,k}^q}{} \leq c \int^{+ \infty}_0 \left(t^{-\theta}\mathcal{K}_{p,k}(t,f; \beta_0,\infty ; \beta_1,\infty)
\right)^q  \frac{dt}{t}<+\infty,
\end{gather*}which proves the
result.
 When $q=+\infty$, we make the usual modif\/ication.

 For $1\leq s\leq q_0, q_1$ Remark~\ref{remark3.2} gives
\begin{gather*}
(\mathcal{B}\mathcal{D}_{p,s}^{\beta_0,k},
\mathcal{B}\mathcal{D}_{p,s}^{\beta_1,k})_{\theta,q} \subset
(\mathcal{B}\mathcal{D}_{p,q_0}^{\beta_0,k}, \mathcal{B}\mathcal{D}_{p,q_1}^{\beta_1,k})_{\theta,q}
 \subset   (\mathcal{B}\mathcal{D}_{p,\infty}^{\beta_0,k},
\mathcal{B}\mathcal{D}_{p,\infty}^{\beta_1,k})_{\theta,q}\subset
\mathcal{B}\mathcal{D}_{p,q}^{\beta,k}.
\end{gather*}Then in order to complete the proof of the theorem we have to show only that
\[
\mathcal{B}\mathcal{D}_{p,q}^{\beta,k}\subset(\mathcal{B}\mathcal{D}_{p,s}^{\beta_0,k},
  \mathcal{B}\mathcal{D}_{p,s}^{\beta_1,k})_{\theta,q}
\quad\mbox{for}\quad 1\leq s\leq q.
\]
 Suppose that $\beta_0>\beta_1$
again. Let $q<+\infty$, we have
\begin{gather*}
\int^{+ \infty}_0
\left(t^{-\theta}\mathcal{K}_{p,k}(t,f; \beta_0,s ; \beta_1,s)
\right)^q  \frac{dt}{t} =   \int^{1}_0  +\int^{+ \infty}_1  = I_1+I_2.
\end{gather*}
Since $\beta>\beta_1$, by Remark~\ref{remark3.2} we get
\[
\mathcal{K}_{p,k}(t,f; \beta_0,s ; \beta_1,s)\leq c t \|f
\|_{\mathcal{ B}\mathcal{D}_{p,s}^{\beta_1,k}}\leq c t \|f
\|_{\mathcal{ B}\mathcal{D}_{p,q}^{\beta,k}},
\]
 hence we deduce
\[
I_1\leq c \|f \|_{\mathcal{ B}\mathcal{D}_{p,q}^{\beta,k}}^q.
\]
To estimate $I_2$ take $f_0= {\sum\limits_{j=0}^{l}\varphi_j
\ast_k f}$ and $f_1= {\sum\limits_{j=l+1}^{+\infty}\varphi_j
\ast_k f}$. Using the properties of the sequence $(\varphi_j)_{j\in
\mathbb{N}}$ we obtain
\[
\|f_0\|_{\mathcal{ B}\mathcal{D}_{p,s}^{\beta_0,k}}^s\leq c  \sum_{j=0}^{l+1}2^{j\beta_0s}\|\varphi_j \ast_k f\|_{p,k}^s\qquad \mbox{and}\qquad  \|f_1
\|_{\mathcal{ B}\mathcal{D}_{p,s}^{\beta_1,k}}^s\leq c
\sum_{j=l}^{+\infty}2^{j\beta_1s}\|\varphi_j \ast_k f\|_{p,k}^s.
\]
Hence we can write
\begin{gather*}
I_2  \leq  c \sum _{l=0}^{+\infty}2^{-\theta q
l(\beta_0-\beta_1)}\Big(\mathcal{K}_{p,k}(2^{l(\beta_0-\beta_1)}
,f; \beta_0,s ; \beta_1,s)\Big)^q \\
\phantom{I_2}{} \leq   c\sum _{l=0}^{+\infty}2^{-\theta q
l(\beta_0-\beta_1)} \!\left[\bigg(\sum_{j=0}^{l+1}2^{j\beta_0s}\|\varphi_j
\ast_k f\|_{p,k}^s\bigg)^{1/s}
+ 2^{l(\beta_0-\beta_1)}\bigg(\sum_{j=l}^{+\infty}2^{j\beta_1s}\|\varphi_j \ast_k f\|_{p,k}^s\bigg)^{1/s}\right]^{q}\!\! \\
\phantom{I_2}{} \leq   c \sum _{l=0}^{+\infty}2^{q
l\beta}\left[\sum_{j=0}^{l+1}2^{(j-l)\beta_0s}\|\varphi_j \ast_k
f\|_{p,k}^s + \sum_{j=l}^{+\infty}2^{(j-l)\beta_1s}\|\varphi_j
\ast_k f\|_{p,k}^s\right]^{q/s}.
\end{gather*}

For $s=q$ it is easy to see that  $I_2 \leq c  \|f
\|_{\mathcal{ B}\mathcal{D}_{p,q}^{\beta,k}}^q.$

For $ s<q$ we take $u>s$ such that
$ {\frac{s}{q}+\frac{s}{u}=1}$ and
$\beta_1<\alpha_1<\beta <\alpha_0<\beta_0$, then by H\"older's
inequality we have
\begin{gather*} I_2  \leq  c  \sum
_{l=0}^{+\infty}2^{ q
l(\beta-\beta_0)}\bigg(\sum_{j=0}^{l+1}2^{(\beta_0-\alpha_0)ju}\bigg)^{q/u}
\bigg(\sum_{j=0}^{l+1}2^{\alpha_0jq}\|\varphi_j \ast_k f\|_{p,k}^q\bigg)\\
\phantom{I_2  \leq}{}  + c \sum _{l=0}^{+\infty}2^{ q
l(\beta-\beta_1)}\bigg(\sum_{j=l}^{+\infty}2^{(\beta_1-\alpha_1)ju}\bigg)^{q/u}
\bigg(\sum_{j=l}^{+\infty}2^{\alpha_1jq}\|\varphi_j \ast_k f\|_{p,k}^q\bigg)\\
\phantom{I_2}{} \leq   c \sum _{l=0}^{+\infty}2^{ q l(\beta-\alpha_0)}
\sum_{j=0}^{l+1}2^{\alpha_0jq}\|\varphi_j \ast_k f\|_{p,k}^q + c
\sum _{l=0}^{+\infty}2^{ q l(\beta-\alpha_1)}
\sum_{j=l}^{+\infty}2^{\alpha_1jq}\|\varphi_j \ast_k f\|_{p,k}^q\\
\phantom{I_2}{} \leq   c \sum_{j=0}^{+\infty}2^{\alpha_0jq}\|\varphi_j \ast_k
f\|_{p,k}^q\sum _{l=j-1}^{+\infty}2^{ q l(\beta-\alpha_0)}
 + c \sum_{j=0}^{+\infty}2^{\alpha_1jq}\|\varphi_j \ast_k f\|_{p,k}^q\sum _{l=0}^{j}2^{ q l(\beta-\alpha_1)}\\
\phantom{I_2}{} \leq   c \|f \|_{\mathcal{ B}\mathcal{D}_{p,q}^{\beta,k}}^q.
\end{gather*}
Finally we deduce
\[
\int^{+ \infty}_0 \left(t^{-\theta}\mathcal{K}_{p,k}(t,f; \beta_0,s ; \beta_1,s)
\right)^q  \frac{dt}{t}\leq  c \|f \|_{\mathcal{
B}\mathcal{D}_{p,q}^{\beta,k}}^q.
\]
 Here when $q=+\infty$ we make the usual modif\/ication.
Our theorem is proved.
\end{proof}

\begin{theorem}\label{theorem3.2} Let  $ \beta>0$ and $1\leq p , q \leq + \infty$.
Then for all $\phi \in \mathcal{A},$ we have
\[
\mathcal{B}\mathcal{D}_{p,q}^{\beta,k} \subset \mathcal{C}_{
p,q}^{\phi ,\beta,k}.
\]
\end{theorem}

\begin{proof} For $\phi \in \mathcal{A}$ and
$1\leq u\leq 2$, we get  $\mbox{supp}\,\mathcal{F}_k(\phi_{2^{-j}u})\subset
A_j$, $\forall\,j\in \mathbb{Z}.$ Then we can write
$\mathcal{F}_k(\phi_{2^{-j}u})=\mathcal{F}_k(\phi_{2^{-j}u})(
\mathcal{F}_k(\varphi_{j-1})+\mathcal{F}_k(\varphi_{j})+\mathcal{F}_k(\varphi_{j+1}))$,
which gives $\phi_{2^{-j}u}=\phi_{2^{-j}u}\ast_k(
\varphi_{j-1}+\varphi_{j}+\varphi_{j+1})$, $\forall\,j\in
\mathbb{Z}.$

 Let $f\in \mathcal{B}\mathcal{D}_{p,q}^{\beta,k}$ for
$1\leq  q < + \infty$, we can assert that
\begin{gather*}
\int^{+ \infty}_0 \left(\frac{\|f \ast_k \phi_t
\|_{p,k}}{t^\beta}\right)^q   \frac{dt}{t}
 \leq {\sum_{j\in \mathbb{Z}}\int^2_1\left(\frac{\|f
\ast_k \phi_{2^{-j}u} \|_{p,k}}{(2^{-j}u)^\beta}\right)^q
\frac{du}{u}}.
\end{gather*}
Using  H\"older's inequality for $j\in \mathbb{Z}$ we get
\[
{\|f \ast_k \phi_{2^{-j}u} \|_{p,k}^q}\leq \| \phi  \|_{1,k}^q
3^{q-1} {\sum_{s=j-1}^{j+1}}\|\varphi_s\ast_k
f\|_{p,k}^q,
\]
hence we obtain
\begin{gather*}
\int^{+ \infty}_0 \left(\frac{\|f \ast_k \phi_t
\|_{p,k}}{t^\beta}\right)^q   \frac{dt}{t}  \leq  c \| \phi
\|_{1,k}^q {\sum_{s\in
\mathbb{Z}}\int^2_1\left(\frac{\|\varphi_{s} \ast_k f
\|_{p,k}}{(2^{-s}u)^\beta}\right)^q   \frac{du}{u}}\\
 \phantom{\int^{+ \infty}_0 \left(\frac{\|f \ast_k \phi_t
\|_{p,k}}{t^\beta}\right)^q   \frac{dt}{t}}{}
 \leq  c \sum_{s\in\mathbb{Z}}(2^{s\beta}\|\varphi_{s}\ast_k
f\|_{p,k})^q  <+\infty.
\end{gather*}
Here when $q=+\infty$, we make the usual modif\/ication. This
completes the proof.
\end{proof}

\begin{theorem}\label{theorem3.3} Let  $
\beta>0 $ and $1\leq p  , q \leq + \infty$, then for $\phi \in
\mathcal{A}$
 such that $ {\sum\limits_{j\in
\mathbb{Z}}\mathcal{F}_k(\phi_{2^{-j}u})(x) =1}$, for all $1\leq
u\leq 2$ and $x\in \mathbb{R}^d$ we have
\[
 \mathcal{C}_{
p,q}^{\phi ,\beta,k}= \mathcal{B}\mathcal{D}_{p,q}^{\beta,k}.
\]
\end{theorem}

\begin{proof} By Theorem~\ref{theorem3.2} we have only to show that $
\mathcal{C}_{ p,q}^{\phi ,\beta,k}\subset
\mathcal{B}\mathcal{D}_{p,q}^{\beta,k} $. Let $\phi \in \mathcal{A}$
such that $ {\sum\limits_{j\in
\mathbb{Z}}\mathcal{F}_k(\phi_{2^{-j}u})(x) =1}$,  for  $x\in
\mathbb{R}^d$ and $1\leq u\leq 2$. Then we can assert that
\[
\mathcal{F}_k(\varphi_{j})=\mathcal{F}_k(\varphi_{j})(\mathcal{F}_k(\phi_{2^{-j-1}u})+
\mathcal{F}_k(\phi_{2^{-j}u})+\mathcal{F}_k(\phi_{2^{-j+1}u})),
\]
this implies that $\varphi_{j}=\varphi_{j}\ast_k(\phi_{2^{-j-1}u}+
\phi_{2^{-j}u}+\phi_{2^{-j+1}u})$, $\forall\,j\in \mathbb{Z}.$

Let $f\in \mathcal{C}_{ p,q}^{\phi ,\beta,k}$ for $1\leq  q < +
\infty$, using  H\"older's inequality for $j\in \mathbb{Z}$ and the
property $ii)$ of the sequence of functions $ (\varphi_{j})_{j\in
\mathbb{Z}}$ we get
\begin{gather*} {\|\varphi_{{j}}\ast_k f   \|_{p,k}^q} \leq  \| \varphi _j \|_{1,k}^q
3^{q-1} {\sum_{s=j-1}^{j+1}}\| f\ast_k
\phi_{2^{-s}u}\|_{p,k}^q \leq c
 {\sum_{s=j-1}^{j+1}}\| f\ast_k \phi_{2^{-s}u}\|_{p,k}^q
. \end{gather*} Integrating with respect to $u$ over $(1, 2)$ we
obtain
\begin{gather*}
\big(2^{j\beta}\|\varphi_{j }\ast_k f\|_{p,k}\big)^q   \leq
 c {\sum_{s=j-1}^{j+1}}\int^2_1\left(\frac{\|f \ast_k
\phi_{2^{-s}u} \|_{p,k}}{(2^{-s}u)^\beta}\right)^q   \frac{du}{u}.
\end{gather*}
Hence
\[
\sum_{j\in\mathbb{Z}}\big(2^{j\beta}\|\varphi_{j }\ast_k
f\|_{p,k}\big)^q \leq c \int^{+ \infty}_0 \left(\frac{\|f \ast_k \phi_t
\|_{p,k}}{t^\beta}\right)^q   \frac{dt}{t} <+\infty.
\] When $q=+\infty$ we make the usual modif\/ication. Our result is proved.
\end{proof}

\begin{remark}\label{remark3.3} We observe that the spaces $ \mathcal{C}_{
p,q}^{\phi ,\beta,k}$ are independent of the specif\/ic selection of
$\phi\in \mathcal{A}$ satisfying the assumption of Theorem~\ref{theorem3.3}.
\end{remark}

\begin{remark}\label{remark3.4} In the case $d=1$, $W= \mathbb{Z}_2$,
$\alpha>-\frac{1}{2}$ and
\[
T_1(f)(x)= \frac{df}{dx}(x) + \frac{2\alpha+1}{x} \left[\frac{f(x)
- f(-x)}{2}\right],\qquad f \in \mathcal{E}( \mathbb{R}),
\] we can
characterize the Besov--Dunkl spaces by dif\/ferences using the Dunkl
translation operators. Observe that
\[
\bigg\{\phi \in \mathcal{A}:
\displaystyle{\sum_{j\in \mathbb{Z}}\mathcal{F}_k(\phi_{2^{-j}u})(x)
=1},\, \forall\;1\leq u\leq 2,\,\forall\,x\in \mathbb{R}
\bigg\}\subset \mathcal{H},
\]
 where $ {\mathcal{H}=\Big\{
\phi \in
\mathcal{S}_\ast(\mathbb{R}) : \int_0^{+\infty}\phi(x)d\mu_\alpha(x)=0
\Big\}}$ with $d\mu_\alpha(x) =
 {\frac{|x|^{2\alpha+1}}{2^{\alpha +1}\Gamma(\alpha
+1)}\, dx}$ and $\mathcal{S}_\ast(\mathbb{R})$ the space of even
Schwartz functions on $\mathbb{R}$. Then we can assert from Theorem~\ref{theorem3.3}
 and~\cite[Theorem~3.6]{[4]} that for $1< p<+\infty$, $1\leq
q\leq+\infty$ and $0<\beta<1$ we have
\[
\mathcal{B}\mathcal{D}_{p,q}^{\beta,k}
=BD^{p,q}_{\alpha,\beta}\subset   \widetilde{
 B}D^{p,q}_{\alpha,\beta},
 \]
where $BD^{p,q}_{\alpha,\beta}$ is the subspace of functions $f \in
L^p( \mu_\alpha)$ satisfying
\[
\left(\int^{+ \infty}_0
\left(\frac{w_{p,\alpha}(f)(x)}{x^\beta}\right)^q
\frac{dx}{x}\right)^{\frac{1}{q}} < + \infty \qquad
\mbox{if}\quad q < +\infty
\] and
\[
 \sup_{x\in (0,+\infty)}\frac{w_{p,\alpha}(f)(x)}{x^\beta} <
+\infty \qquad \mbox{if} \quad q=+\infty ,
\]
 with
$w_{p,\alpha}(f)(x)= \|\tau_x(f) + \tau_{-x}(f) - 2f\|_{p,\alpha}$.
For the space $\widetilde{B}D^{p,q}_{\alpha,\beta}$ we replace
$w_{p,\alpha}(f)(x)$ by $\widetilde{w}_{p,\alpha}(f)(x)=
\|\tau_x(f)- f\|_{p,\alpha}$.

 Note that when $f$ is an even
function in $L^p( \mu_\alpha)$ we have
$\tau_x(f)(y)=\tau_{-x}(f)(-y)$ for $x,y\in\mathbb{R}$, then we get
\[
f\in BD^{p,q}_{\alpha,\beta}\Longleftrightarrow f\in \widetilde{
 B}D^{p,q}_{\alpha,\beta}.
 \]
 \end{remark}

\section[Integrability of the Dunkl transform of function in Besov-Dunkl space]{Integrability of the Dunkl transform of function\\ in Besov--Dunkl space}\label{section4}

$ $ In this section, we establish further results concerning
integrability of the Dunkl transform of function $f$ on
$\mathbb{R}^d$, when $f$ is in a suitable Besov--Dunkl space.

In the following lemma we prove the Hardy--Littlewood inequality for
the Dunkl transform.

\begin{lemma}\label{lemma4.1}
If $f \in L_k^p(
\mathbb{R}^d)$ for some $1 < p \leq 2$, then
\begin{gather}\label{eq11}
\int_{\mathbb{R}^d} \|x\|^{2(\gamma + \frac{d}{2})(p-2)}
|\mathcal{F}_k(f)(x)|^p w_k(x)dx \leq c
\|f\|_{p,k}^p.\end{gather}
\end{lemma}

\begin{proof} To see~\eqref{eq11} we will make
use of the Marcinkiewicz interpolation theorem (see~\cite{[20]}). For $f
\in L_k^p( \mathbb{R}^d)$ with $1 \leq p \leq 2$ consider the
operator
\[
\mathcal{L}(f) (x) = \|x\|^{2(\gamma + \frac{d}{2})} \mathcal{F}_k(f)(x),\quad x \in \mathbb{R}^d.
\]
For every $f \in L^2( \mathbb{R}^d)$ we have from Plancherel's
theorem
\begin{gather}\label{eq12}
\left(\int_{\mathbb{R}^d}|\mathcal{L}(f)(x)|^2 \frac{w_k(x)}
{\|x\|^{4(\gamma +\frac{d}{2})}}dx\right)^{1/2} =
\|\mathcal{F}_k(f)\|_{2,k}= c^{-1}_k \|f\|_{2,k},
\end{gather}
Moreover, according to \eqref{eq1} and \eqref{eq2} we get for $\lambda \in ]0, +
\infty)$ and $f \in L^1_k(\mathbb{R}^d)$
\begin{gather}
\int_{\{x \in \mathbb{R}^d:\mathcal{L}(f)(x)|
> \lambda\}} \frac{w_k(x)} {\|x\|^{4(\gamma + \frac{d}{2})}} dx \leq  \int_{\|x\| >
(\frac{\lambda}{\|f\|_{1,k}})^{{\frac{1}{2(\gamma +
\frac{d}{2})}}}} \frac{w_k(x)}{\|x\|^{4(\gamma + \frac{d}{2})}}dx\nonumber\\
\phantom{\int_{\{x \in \mathbb{R}^d:\mathcal{L}(f)(x)| > \lambda\}} \frac{w_k(x)} {\|x\|^{4(\gamma + \frac{d}{2})}} dx}{}  \leq  c\int^{+ \infty}_{(\frac{\lambda}
{\|f\|_{1,k}})^{{\frac{1}{2(\gamma+\frac{d}{2})}}}}\frac{r^{2\gamma+d-1}}
{r^{4(\gamma + \frac{d}{2})}}dr \leq
c\,\frac{\|f\|_{1,k}}{\lambda} .\label{eq13}
\end{gather} Hence by \eqref{eq12} and
\eqref{eq13} $\mathcal{L}$ is an operator of strong-type $(2,2)$ and
weak-type $(1,1)$ between the spaces $( \mathbb{R}^d, w_k(x)dx)$ and
$( \mathbb{R}^d, \frac{w_k(x)} {\|x\|^{4(\gamma +
\frac{d}{2})}}dx)$.

Using Marcinkiewicz interpolation's theorem we can assert
that $\mathcal{L}$ is an operator of strong-type $(p,p)$ for $1 < p
\leq 2$, between the spaces under consideration. We conclude that
\begin{gather*}
\int_{ \mathbb{R}^d} |\mathcal{L}(f)(x)|^p
\frac{w_k(x)}{\|x\|^{4(\gamma + \frac{d}{2}})}dx  =
\int_{\mathbb{R}^d} \|x\|^{2(\gamma +
\frac{d}{2})(p-2)} |\mathcal{F}_k(f)(x)|^p w_k(x)dx
 \leq  c  \|f\|_{p,k}^p ,
\end{gather*}
thus we obtain the result.
\end{proof}

Now in order to prove the following two theorems we denote by
$\mathcal{\widetilde{A}}$ the subset of functions $\phi$ in
$\mathcal{A}$ such that
\begin{gather}\label{eq14}
\exists\,c>0 ;\quad |\mathcal{F}_k(\phi)(x)|\geq c \|x\|^2\qquad \mbox{if}\quad 1\leq \|x\|\leq 2 .
\end{gather}
Let $ \beta
 >0 $ and $1\leq p  , q \leq + \infty$. From Theorem~\ref{theorem3.2} we
have obviously for all $\phi \in \mathcal{\widetilde{A}}$,
\begin{gather}\label{eq15}
\mathcal{B}\mathcal{D}_{p,q}^{\beta,k} \subset \mathcal{C}_{
p,q}^{\phi ,\beta,k}.
\end{gather} For $1\leq p \leq 2$ we take
$p'$ such that $\frac{1}{p}+\frac{1}{p'}=1$. We recall that
$\mathcal{F}_k(f)\in L^{p'}_{k}(\mathbb{R}^d)$ for all $f \in
L^p_k(\mathbb{R}^d)$.

\begin{theorem}\label{theorem4.1}
Let $1<p \leq 2$. If
$f\in \mathcal{B}\mathcal{D}_{p,1}^{ \frac{
2(\gamma+\frac{d}{2})}{p},k}$, then
\[
\mathcal{F}_k(f) \in L^1_k( \mathbb{R}^d).
\]
\end{theorem}

\begin{proof} Let $f \in \mathcal{B}\mathcal{D}_{p,1}^{ \frac{
2(\gamma+\frac{d}{2})}{p},k}$ with $1<p\leq 2$. For $ \phi \in
\mathcal{\widetilde{A}}$ we can write from~\eqref{eq6} and for $t \in
(0,+\infty)$, $\mathcal{F}_k( f\ast_k \phi_t)(x) =
\mathcal{F}_k(f)(x) \mathcal{F}_k(\phi_t)(x)$,  a.e.\   $x
\in \mathbb{R}^d.$  From Lemma~\ref{lemma4.1}
  we obtain
  \[
  \int_{\mathbb{R}^d}
|\mathcal{F}_k(f)(x)|^p  |\mathcal{F}_k(\phi_t)(x)|^p
 \|x\|^{2(\gamma + \frac{d}{2})(p-2)}w_k(x)dx  \leq c \|f\ast_k
\phi_t\|_{p,k}^p  .
\] By \eqref{eq14} we get $
|\mathcal{F}_k(\phi_t)(x)|\geq c \|tx\|^2$ if $1\leq
\|tx\|\leq 2 $, then we can assert that
 \begin{gather}\label{eq16}
 t ^{2} \left(\int_{\frac{1}{t}\leq\|x\|\leq \frac{2}{t}} |\mathcal{F}_k(f)(x)|^p
 \|x\|^{2(\gamma + \frac{d}{2})(p-2)+2p}w_k(x)dx\right)^{1/p}
 \leq c \|f\ast_k \phi_t\|_{p,k}.
 \end{gather}
 Then by H\"older's inequality, \eqref{eq1} and \eqref{eq16} we have
\begin{gather*}
\int_{{\frac{1}{t}\leq\|x\|\leq
\frac{2}{t}}}\|x\| |\mathcal{F}_k(f)(x)|  w_k(x)dx \leq  c
 \frac{\|f\ast_k \phi_t\|_{p,k}}{t^{2}}
\left(\int^{\frac{2}{t}}_{\frac{1}{t}} r^{(\frac{1}{1-p})[\,2(\gamma
+ \frac{d}{2})(p-2)+p\, ]}\, r^{2\gamma+d-1}dr\right)^{\frac{1}{p'}}
\\
\phantom{\int_{{\frac{1}{t}\leq\|x\|\leq
\frac{2}{t}}}\|x\| |\mathcal{F}_k(f)(x)|  w_k(x)dx}{} \leq  c   \frac{\|f\ast_k \phi_t\|_{p,k}}{t^{
\frac{2(\gamma+\frac{d}{2})}{p}}}   \frac{1}{t} .
\end{gather*} Integrating with respect to $t$
over $\mathbb{R}_+$, applying Fubini's theorem and using~\eqref{eq15}, it
yields
\begin{eqnarray*}\int_{\mathbb{R}^d}|\mathcal{F}_k(f)(x)| \,w_k(x)dx \leq c \int^{+ \infty}_0
 \frac{\|f\ast_k \phi_t\|_{p,k}}{t^{\frac{2(\gamma+\frac{d}{2})}{p}}} \, \frac{dt}{t} <
+\infty \,.\end{eqnarray*} This
 complete the proof of the theorem.
\end{proof}

\begin{theorem}\label{theorem4.2} Let $\beta >0$ and $1\leq p \leq 2$. If $f\in
\mathcal{B}\mathcal{D}_{p,\infty}^{\beta,k}$, then
\begin{itemize}\itemsep=0pt
\item[$i)$] for $p\neq 1 $ and $ 0<\beta\leq \frac{2(\gamma+\frac{d}{2})}{p}$ we have
\[
\mathcal{F}_k(f) \in L^{s}_{k}(\mathbb{R}^d)\qquad\mbox{provided that}\quad \frac{2(\gamma +\frac{d}{2})p}{\beta
p+2(\gamma+\frac{d}{2})(p-1)} < s \leq p';
\]
\item[$ii)$] for   $ \beta>\frac{2(\gamma+\frac{d}{2})}{p}$ we have  $\mathcal{F}_k(f) \in L^{1}_{k}(\mathbb{R}^d).$
\end{itemize}
\end{theorem}

\begin{proof}
 Let $f\in
\mathcal{B}\mathcal{D}_{p,\infty}^{\beta,k}\;$  with
$1\leq p \leq 2$ and $ \phi \in \mathcal{\widetilde{A}}$.

$i)$ Suppose that $p\neq 1$ and $0<\beta\leq
\frac{2(\gamma+\frac{d}{2})}{p}$. Using \eqref{eq3} and \eqref{eq6} we have for
$t\in (0,+\infty )$
\[
\|\mathcal{F}_k(f\ast_k \phi_t) \|_{p',k}=
\|\mathcal{F}_k(f)\mathcal{F}_k(\phi_t) \|_{p',k}\leq c \|f\ast_k
\phi_t \|_{p,k} .
\] Then from \eqref{eq14} and \eqref{eq15} we obtain
 \begin{gather}\label{eq17}
 t^2\left(\int_{\frac{1}{t}\leq \|x\|\leq \frac{2}{t}} |\mathcal{F}_k(f)(x)|^{p'}  \|x\|^{2p'} w_k(x)dx\right)^{1/p'}
 \leq c  \|f\ast_k \phi_t \|_{p,k}\leq  c t^\beta .
 \end{gather}
Let $ s\in \big] {\frac{2(\gamma +\frac{d}{2})p}{\beta
p+2(\gamma+\frac{d}{2})(p-1)}}  ,  p'  \big]$. Since
$\mathcal{F}_k(f)\in L^{p'}_{k}(\mathbb{R}^d)$, we have only to show
the case $s\neq p'$. For $t\geq1$ put $G_t$ the set of $x$ in
$\mathbb{R}^d $ such that $ \frac{1}{t^{1/s}}\leq\|x\|\leq
\frac{2}{t^{1/s}}$.
By H\"older's inequality, \eqref{eq1} and~\eqref{eq17} we have
\begin{gather*}
\int_{G_t}|\mathcal{F}_k(f)(x)|^{s}\, \|x\|^{s}
 w_k(x)dx \\
\qquad{} \leq
\left(\int_{G_t}|\mathcal{F}_k(f)(x)|^{p'} \|x\|^{2p'}w_k(x)dx\right)^{s/p'}
\left(\int_{G_t}\|x\|^{\frac{-p's}{p'-s}}w_k(x)dx\right)^{ 1-
\frac{s}{p'}}\\
\qquad{} \leq  c t^{\beta-2}
\left(\int_{\frac{1}{t^{1/s}}}^{\frac{2}{t^{1/s}}}
r^{2\gamma+d-1-\frac{p's}{p'-s}}dr\right)^{ 1- \frac{s}{p'}}  \leq c
t^{-1+\beta - 2(\gamma+\frac{d}{2})(\frac{1}{s} - \frac{1}{p'})}.
\end{gather*}
Integrating with respect to $t$ over $(0,1)$
and applying Fubini's theorem, it yields
\[
\int_{\|x\|\geq 1} |\mathcal{F}_k(f)(x)|^s w_k(x)dx \leq c\int^{1}_0 t^{-1+\beta - 2(\gamma+\frac{d}{2})(\frac{1}{s} -
\frac{1}{p'})}dt <+\infty .
\] Since $L^{p'}_k(B(0,1),w_k(x)dx)
\subset L^{s}_k(B(0,1),w_k(x)dx)$ we deduce that $\mathcal{F}_k(f)$
is in $L^{s}_k(\mathbb{R}^d)$.

$ii)$ Assume now
$\beta>\frac{2(\gamma+\frac{d}{2})}{p}$.
 For $p\neq1$
by proceeding in the same manner as in the proof of $i)$ with $s=1$,
we obtain the desired result.

For $p=1$, using \eqref{eq3} and \eqref{eq6}, we have for $t\in (0,+\infty )$
\[
\|\mathcal{F}_k(f\ast_k \phi_t) \|_{\infty,k}=
\|\mathcal{F}_k(f)\mathcal{F}_k(\phi_t) \|_{\infty,k}\leq
c \|f\ast_k \phi_t \|_{1,k} .
\] Then from \eqref{eq14} and \eqref{eq15} we obtain
\begin{gather}\label{eq18}
t^2\|h_t \mathcal{F}_k(f) \|_{\infty,k}
 \leq c  \|f\ast_k \phi_t \|_{1,k}\leq  c t^\beta,
  \end{gather}
  where $h_t(x)=\chi_t(x)\|x\|^2$
 with $\chi_t$ is the characteristic function of the set $\{x\in \mathbb{R}^d: {\frac{1}{t}\leq\|x\|\leq
\frac{2}{t}}\}$.

 By H\"older's inequality, \eqref{eq1} and \eqref{eq18} we have
\begin{gather*}\int_{{\frac{1}{t}\leq\|x\|\leq
\frac{2}{t}}}|\mathcal{F}_k(f)(x)|  \|x\| w_k(x)dx  \leq
\|h_t \mathcal{F}_k(f) \|_{\infty,k}\int_{\mathbb{R}^d}
|\chi_t(x)| \|x\|^{-1} w_k(x)dx\\
\phantom{\int_{{\frac{1}{t}\leq\|x\|\leq
\frac{2}{t}}}|\mathcal{F}_k(f)(x)|  \|x\| w_k(x)dx }{}
\leq  c t^{\beta-2}
\int_{\frac{1}{t}}^{\frac{2}{t}} r^{2\gamma+d-2}\, dr \leq  c  t^{
\beta - 2(\gamma+\frac{d}{2})-1}.
\end{gather*} Integrating with respect to $t$ over $(0,1)$
and applying Fubini's theorem we obtain
\[
\int_{\|x\|\geq 1}
|\mathcal{F}_k(f)(x)| w_k(x)dx \leq c\int^{1}_0 t^{ \beta -
2(\gamma+\frac{d}{2})-1}dt <+\infty .
\] Since
$L^{\infty}_k(B(0,1),w_k(x)dx) \subset L^{1}_k(B(0,1),w_k(x)dx)$ we
deduce that $\mathcal{F}_k(f)$ is in $L^{1}_k(\mathbb{R}^d)$. Our
theorem is proved.
\end{proof}

\begin{remark}\label{remark4.1} \quad {}
\begin{enumerate}\itemsep=0pt \item [1.] For $\beta>0$, $1\leq p\leq 2$ et $1\leq q\leq +\infty$, using Remark~\ref{remark3.2}, the results
of Theorem~\ref{theorem4.2} are true for
$\mathcal{B}\mathcal{D}_{p,q}^{\beta,k}$.

\item [2.] From Remark~\ref{remark3.2}
we get $
\mathcal{B}\mathcal{D}_{p,\infty}^{\beta,k}\subset\mathcal{B}\mathcal{D}_{p,1}^{
\frac{ 2(\gamma+\frac{d}{2})}{p},k}$ for $\beta>
\frac{2(\gamma+\frac{d}{2})}{p}$. Using Theorem~\ref{theorem4.1} we recover the
result of Theorem~\ref{theorem4.2},~$ii)$ with $1<p\leq2$.

\item [3.] Let $\beta>
2(\gamma+\frac{d}{2})$, by Theorem~\ref{theorem4.2},~$ii)$ we can assert that
\begin{enumerate}\itemsep=0pt
\item [$i)$] $  \mathcal{B}\mathcal{D}_{1,\infty}^{\beta,k} $ is an
example of space where we can apply the inversion formula;

\item [$ii)$] $\mathcal{B}\mathcal{D}_{1,\infty}^{\beta,k}$ is contained
in $L^{1}_k(\mathbb{R}^d)\cap L^{\infty}_k(\mathbb{R}^d)$ and hence
is a subspace of $L^{2}_k(\mathbb{R}^d)$. By~\eqref{eq4} we obtain for
$f\in \mathcal{B}\mathcal{D}_{1,\infty}^{\beta,k} $
\[
\tau_y(f)(x) =c_k\int_{\mathbb{R}^d}\mathcal{F}_k(f)(\xi) E_k(ix, \xi)E_k(-iy,
\xi)w_k(\xi)d\xi,\qquad  x ,y\in \mathbb{R}^d.
\]
\end{enumerate}
\end{enumerate}
\end{remark}

\subsection*{Acknowledgements}
The authors thank the referees for their remarks and suggestions.
Work supported by the DGRST research project 04/UR/15-02 and the
program CMCU 07G 1501.

\pdfbookmark[1]{References}{ref}
\LastPageEnding


\begin{thebibliography}{99}

\footnotesize\itemsep=0pt

\bibitem{[1]} Abdelkef\/i C., Sif\/i M., On the uniform convergence of
partial Dunkl integrals in Besov--Dunkl spaces, {\it  Fract. Calc. Appl. Anal.} {\bf 9} (2006), 43--56.

\bibitem{[2]}  Abdelkef\/i C., Sif\/i M., Further results of integrability for the Dunkl
transform, {\it Commun. Math. Anal.} {\bf 2}
(2007), 29--36.

\bibitem{[3]} Abdelkef\/i C., Dunkl transform on Besov spaces
and Herz spaces, {\it  Commun. Math. Anal.} {\bf 2}
(2007), 35--41.

\bibitem{[4]}  Abdelkef\/i  C., Sif\/i M., Characterization of
Besov spaces for the Dunkl operator on the real line, {\it  JIPAM. J. Inequal. Pure Appl. Math.} {\bf 8} (2007), no.~3, Article~73, 11~pages.

\bibitem{[5]}  Bergh  J., L\"{o}fstr\"{o}m J.,
Interpolation spaces. An introduction, {\it Grundlehren der Mathematischen Wissenschaften}, no.~223, Springer-Verlag, Berlin~-- New York,  1976.

\bibitem{[6]} Besov O.V., On a family of function
spaces in connection with embeddings and extentions, {\it Trudy. Mat. Inst. Steklov.} {\bf  60} (1961), 42--81 (in Russian).

\bibitem{[7]} Betancor  J.J.,  Rodr\'{\i}guez-Mesa L., Lipschitz--Hankel spaces and partial Hankel
integrals, {\it Integral Transform. Spec. Funct.} {\bf  7}  (1998), 1--12.

\bibitem{[8]}  de Jeu  M.F.E., The Dunkl transform, {\it Invent. Math.} {\bf 113} (1993),
147--162.

\bibitem{[9]} Dunkl C.F., Dif\/ferential-dif\/ference operators
associated to ref\/lection groups, {\it Trans. Amer. Math. Soc.} {\bf  311}
(1989), 167--183.

\bibitem{[10]} Dunkl C.F., Integral kernels with ref\/lection
group invariance, {\it Canad. J. Math.} {\bf 43}, (1991), 1213--1227.

\bibitem{[11]}
 Dunkl C.F., Hankel transforms associated to f\/inite ref\/lection
groups, in  Proc. of Special Session on Hypergeometric Functions on
Domains of Positivity, Jack Polynomials and Applications (Tampa, 1991), {\it Contemp. Math.} {\bf 138} (1992),
123--138.

\bibitem{[12]} Dunkl  C.F., Xu Y., Orthogonal polynomials of
several variables, {\it Encyclopedia of Mathematics and its Applications}, Vol.~81, Cambridge University Press, Cambridge, 2001.

\bibitem{[13]} Giang  D.V., M\'{o}rciz F., On the uniform and the absolute
convergence of Dirichlet integrals of functions in Besov-spaces,
{\it Acta Sci. Math. (Szeged)} {\bf 59} (1994), 257--265.

\bibitem{[14]} Kamoun L.,
Besov-type spaces for the Dunkl operator on the real line, {\it J. Comput. Appl. Math.} {\bf  199} (2007), 56--67.


\bibitem{[15]}  Pelczynski  A., Wojciechowski M., Molecular decompositions
and embbeding theorems for vector-valued Sobolev spaces with
gradient norm, {\it Studia Math.} {\bf 107} (1993), 61--100.

\bibitem{[16]}  R\"osler M., Generalized Hermite polynomials and the heat equation for
Dunkl operators, {\it Comm. Math. Phys.} {\bf 192} (1998), 519--542, \href{http://arxiv.org/abs/q-alg/9703006}{q-alg/9703006}.

\bibitem{[17]} R\"osler M., Positivity of Dunkl's intertwining operator, {\it Duke Math. J.} {\bf 98} (1999), 445--463, \href{http://arxiv.org/abs/q-alg/9710029}{q-alg/9710029}.

\bibitem{[18]} R\"osler  M., Voit M., Markov processes with Dunkl operators, {\it Adv. in Appl. Math.} {\bf 21} (1998), 575--643.

\bibitem{[19]} Thangavelyu  S., Xu Y., Convolution operator and
maximal function for Dunkl transform, {\it J. Anal. Math.} {\bf 97} (2005), 25--55.

\bibitem{[20]} Titchmarsh E.C., Introduction to the theory of Fourier
integrals, Clarendon Press, Oxford, 1937.

\bibitem{[21]} Triebel H., Theory
of function spaces, {\it Monographs in Mathematics}, Vol.~78, Birkh\"auser,
Verlag, Basel, 1983.

\bibitem{[22]} Trim\`eche K., The Dunkl intertwining
operator on spaces of functions and distributions and integral
representation of its dual, {\it Integral Transforms Spec. Funct.} {\bf  12}
(2001), 349--374.

\bibitem{[23]} Trim\`eche K., Paley--Wiener theorems for the
Dunkl transform and Dunkl translation operators, {\it Integral Transforms Spec. Funct.} {\bf  13} (2002), 17--38.
\end{thebibliography}
\end{document}